\documentclass[reqno,12pt]{amsart}
\usepackage{amsmath}
\usepackage{amssymb}
\usepackage{cite}
\usepackage{amsmath}
\usepackage{txfonts}
\usepackage{stmaryrd}
\usepackage{amssymb}
\usepackage{amsfonts}
\usepackage{times}
\usepackage{mathrsfs}
\usepackage{amsthm}
\usepackage{enumerate}
\usepackage{comment}
\usepackage{hyperref}
\hypersetup{hidelinks,
 colorlinks=true,
 allcolors=blue,
 pdfstartview=Fit,
 breaklinks=true}
\usepackage{ifthen}
\usepackage[usenames]{color}
\usepackage{color}
\usepackage{CJK}

\newcommand{\D}{{\mathbb D}}

\def\a{\alpha}

\def \psi {u}

\def\msk{\medskip}

\def\bege{\begin{equation}}
\def\ende{\end{equation}}
\def\begr{\begin{eqnarray}}
\def\endr{\end{eqnarray}}
\def\bnum{\begin{enumerate}}
\def\enum{\end{enumerate}}

\theoremstyle{definition}

 \newtheorem{thm}{Theorem}[section]
 
 \newtheorem{lem}{Lemma}[section]
 
 \theoremstyle{definition}
 
 \theoremstyle{remark}

 \numberwithin{equation}{section}

\begin{document}

\title[order boundedness and essential norm of composition operators]
 {order boundedness and essential norm of generalized weighted composition operators on Bergman spaces with doubling weights}

\author[Zuoling Liu]{Zuoling Liu}
\address{Zuoling Liu: Department of Mathematics, Shantou University, Shantou 515063, China;}
\email{16zlliu1@stu.edu.cn}





\date{\today}

\begin{abstract}
In this paper, the order boundedness and essential norm of generalized weighted composition operators on Bergman spaces with doubling weights are characterized. Specially, we estimate the essential norm of these operators on weighted Bergman spaces by using the reduce order method.
\end{abstract}

\maketitle

\section{Introduction}

Let $\mathbb{D} = \{z:\ \lvert z \rvert<1\}$ be the open unit disk in the complex plane $\mathbb{C}$. Let $\omega: \mathbb{D} \rightarrow [0,\infty)$ be an integrable function and radial, that is, $\omega(z)=\omega(|z|)$ for all $z\in \mathbb{D}$. Denote $\hat{\omega}(z)=\int_{|z|}^{1}\omega(s)ds$ for all $z\in \mathbb{D}$. A weight $\omega$ belongs to the class $\hat{\mathcal{D}}$ if $\hat{\omega}(r)\leq C \hat{\omega}(\frac{1+r}{2})$, where $ C=C(\omega)\geq 1$ and $0\leq r < 1$. Furthermore, we write $\omega\in \breve{\mathcal{D}}$ if there exist  constants $\vartheta=\vartheta(\omega)>1$ and $ C=C(\omega)> 1$   such that $\hat{\omega}(r)\geq C \hat{\omega}(1-\frac{1-r}{\vartheta})$ for all $0\leq r < 1$. We denote $\mathcal{D}=\check{\mathcal{D}} \cap \hat{\mathcal{D}}$ and $\omega(E)=\int_{E}\omega dA$ for each measurable set $E\subset \mathbb{D}$.

Let $H(\mathbb{D})$ be the space of analytic functions on $\mathbb{D}$.  For $0<p<\infty$ and the radial weight $\omega$,  the Bergman  space $A^{p}_{\omega}$ associated to $\omega$ is defined by
 \[
 A^{p}_{\omega}= \left\{ f\in H(\mathbb{D}) : \|f\|^p_{A^{p}_{\omega}}=  \int_{\mathbb{D}}|f(z)|^{p} \omega(z)dA(z)  <\infty \right\},
 \]
 where $dA(z)=\frac{1}{\pi}dxdy$ is the normalized Lebesgue area measure on $\mathbb{D}$.  As
usual, let $A^p_\a$ stand  for the classical weighted Bergman space induced by   radial weight $\omega(z)=(1-|z|^2)^\a$, where $-1<\a<\infty$.
$A^p_\omega$  is a Banach space for $1\leq p<\infty$ under the norm $\|\cdot\|_{A^p_\omega}$. 
See \cite{HKZ2000, PR2014} for  the theory of  weighted Bergman spaces.   Let $q>0$ and $\mu$  be a finite positive Borel measure on $\mathbb{D}$. We say that  $f\in L_{\mu}^q$ if the measurable function $f$ satisfies
$$\|f\|^q_{L_{\mu}^q}= \int_{\mathbb{D}}|f(w)|^qd\mu(w) <\infty.  $$


 Suppose $\varphi$ is an analytic map of   $\mathbb{D}$ into itself. Every analytic self-map $\varphi$ induces a composition operator $C_{\varphi}$ on $H(\mathbb{D})$ by $$C_{\varphi}(f)(z)=f(\varphi(z)),~~f\in H(\mathbb{D}),$$
for all $z\in \mathbb{D}$. See \cite{CM1995} and \cite{S1993} for  the theory of  composition operators.

For $n\in \mathbb{N}$, $D^n f=f^{(n)}$ is the differential operator on $H(\mathbb{D})$. Let $n\in  \mathbb{N}\cup\{0\} $ and $u\in H(\mathbb{D})$. The generalized weighted composition operator, denoted by $D^{n}_{\varphi,u}$, is defined by
\[D^{n}_{\varphi,u}(f)(z)=u(z)f^{(n)}(\varphi(z)),~~f\in H(\mathbb{D}).
\]
   The generalized weighted composition operator was  coined by Zhu in \cite{Z2007}. Clearly, if $n=0$ and $u\equiv1$, the operator $D^{n}_{\varphi,u}$ becomes the composition operator $C_{\varphi}$.  When $n=0$,  the operator $D^{n}_{\varphi,u}$ is called the weighted composition operator, usually denoted by $uC_{\varphi}$.  As $n=1$ and $u(z)=\varphi'(z)$, then $D^{n}_{\varphi,u}= DC_{\varphi}$. If $n=1$ and $u(z)\equiv1$, then $D^{n}_{\varphi,u}= C_{\varphi}D$.  The operators $DC_{\varphi}$ and $ C_{\varphi}D$  have been studied in \cite{HP2005,MS2014,MS2015,MS2018}.   For some recent work on generalized weighted composition operators, we refer the interested readers to \cite{Y2015,Z2008,Z2018} and \cite{MG2021}.


Let $\mathbf{X}$ be a quasi-Banach space and $\mu$ be a positive measure on $\D$. Assume $0 < q < \infty$ and let $T : X  \rightarrow L^q_\mu$ be an operator. We say that $T$ is order bounded if $T$ maps the closed unit ball $B_{\mathbf{X}}$ of $\mathbf{X}$ into an order interval of $L^q_\mu$. In other words, there exists a non-negative element $h\in L^q_\mu$ such that $|Tf|\leq h$ almost everywhere with respect to $\mu$ for all $f\in B_{\mathbf{X}}$. This concept has been studied in several references \cite{U2012,SS2019,S2017}.

The order bounded composition operators on Hardy spaces was introduced by Hunziker and Jarchow in \cite{HJ1991}. Motivated by \cite{HJ1991}, Hibschweiler \cite{H2008} characterized order bounded composition operators acting on standard weighted Bergman spaces. Later, Ueki \cite{U2012} considered the order boundedness of weighted composition operators on standard weighted Bergman spaces. Subsequently,  Wolf \cite{E2012} studied order bounded weighted composition operators acting on Bergman spaces with general weights.
Recently, the order boundedness of weighted composition operators acting between Banach spaces like Hardy spaces, weighted Bergman spaces, weighted Dirichlet spaces and derivative Hardy spaces were discussed (see \cite{SS2019,S2017,GKZ2016,LLW2020}).  Motivated by \cite{U2012,E2012}, we investigate the order boundedness of $D^{n}_{\varphi,u}$ on Bergman spaces with doubling weights.

Let $X$ and $ Y$ be Banach spaces.
The essential norm of linear operator $T:X\rightarrow Y$ is defined as
$$\|T\|_{e,~X\rightarrow Y}=\inf_{K}\|T-K\|_{X\rightarrow Y},$$
where $K$ is any compact operator and $\|\cdot\|_{X\rightarrow Y}$ is the operator norm. It is obvious that  $\|T\|_{e,~X\rightarrow Y}=0$ if and only if $T$ is a compact operator.

The study of the essential norm of composition operators on Hardy spaces and Bergman spaces was dated back to Shapiro \cite{S1978}.
\u{C}u\u{c}kovi\'{c} and Zhao extended Shapiro's \cite{S1978} results to standard weighted Bergman spaces and Hardy spaces in \cite{CZ2004, CZ2007}. After their works, Demazeux \cite{D2011} considered the essential norm of weighted composition operators on Hardy spaces in terms of pullback measure  for $1\leq p,q\leq\infty$.
In light of their work, the authors \cite{DLS2020} investigated the boundedness and essential norm of weighted composition operators on Bergman spaces induced by doubling weights. Based on their work and inspired by the idea from \cite{Z2008}, Liu \cite{L2021} studied the  boundedness and compactness of generalized weighted composition operators $D^{n}_{\varphi,u}$ between different Bergman spaces with doubling weights. See \cite{DLS2020, LRW2021, LP2022, L2021} for more results of composition operators on  Bergman spaces $A^{p}_{\omega}$.
The essential norm of composition operators on Bergman spaces with admissible B\'{e}koll\'{e} weights was studied by \cite{SU2014}. Recently,  Esmaeili and Kellay \cite{KK2022} considered the essential norm of weighted composition operators on weighted Bergman spaces.
Many authors considered the essential norm of composition operators on different weighted Bergman spaces, see \cite{ UL2008, LC2013} and references therein.

Motivated by the idea from \cite{DLS2020,KK2022,SU2014}, we estimate the essential norm of generalized weighted composition operators on Bergman spaces with doubling weights. 

In this paper, we denote constants by $C$ which are positive and may differ from one occurrence to the other.
The notation $a \lesssim b$ means that there is a positive constant $C$
such that $a \leq C b$.
The symbol $a \asymp b$ means that both $a\lesssim b$ and $b\lesssim a$ hold.

\section{Preliminary results}
The pseudo-hyperbolic metric $\rho$ on $\mathbb{D}$ is defined as $$\rho(z,w)=|\varphi_{w}(z)|=\Big|\frac{w-z}{1-\bar{w}z}\Big|,$$  for $z, w\in \mathbb{D}$. For $r\in(0,1)$, the pseudo-hyperbolic disk is defined by $$\Delta(w,r)=\{z\in\mathbb{D},~~\rho(z,w)<r\}.$$
For $z\in \mathbb{D}\backslash\{0\}$,
 $$S(z)=\left\{\xi\in \mathbb{D}:~~|z|\leq|\xi|< 1,~~|\arg\xi-\arg z|<\frac{1-|z|}{2}\right\}$$ is called a Carleson square. We set $S(0)=\mathbb{D}$.

\begin{lem}\cite[Lemma 2.1]{L2021}\label{DZ}
Let $\omega \in \mathcal{D}$, $0<p<\infty$ and $n\in \mathbb{N}\bigcup \{0\}$. If $f\in A_{\omega}^p$, then there exists a constant $C=C(\omega)>0$ such that
$$|f^{(n)}(z)|\leq C \frac{||f||_{A_{\omega}^p}}{ (\omega(S(z)))^{1/p} (1-|z|)^n }$$
for all $z\in \mathbb{D}$.
\end{lem}

\begin{lem}\cite[Proposition 3.1]{L2021}\label{D}
Let $0<p\leq q < \infty$, $\omega \in  \mathcal{D}$ and $n\in \mathbb{N}\bigcup \{0\}$. Let $\mu$ be a positive Borel measure on $\mathbb{D}$. Then there exists $r=r(\omega)\in (0,1)$ such that the following statements hold.
\begin{itemize}
\item[(i)]$D^n: A_{\omega}^{p}  \rightarrow L_{\mu}^{q}$ is bounded if and only if
\begin{equation}\label{i}
\sup_{z\in\mathbb{D}} \frac{ \mu(\Delta(z,r))}{ (\omega(S(z)))^{q/p}(1-|z|)^{nq} }<\infty.
\end{equation}
Moreover,
$$
\|D^{n}\|^{q}_{A_{\omega}^{p}  \rightarrow L_{\mu}^{q}}\asymp\sup_{z\in\mathbb{D}}\frac{\mu(\Delta(z,~r))}{(1-|z|)^{nq}(\omega(S(z)))^{q/p}}.
$$
\item[(ii)]$D^n: A_{\omega}^{p}  \rightarrow L_{\mu}^{q}$ is compact if and only if
\begin{equation}\label{ii}
\lim_{|z|\rightarrow 1^{-}}  \frac{ \mu(\Delta(z,r))}{ (\omega(S(z)))^{q/p}(1-|z|)^{nq} }=0.
\end{equation}
\end{itemize}

\end{lem}
 In light of Lemma \ref{D}, we get the following lemma.

 \begin{lem}\label{LEM1}
Let $0<p\leq q<\infty$, $\omega \in  \mathcal{D}$ and $n\in \mathbb{N}\bigcup \{0\}$. Assume that  $\mu$ is a positive Borel measure on $\mathbb{D}$, $r=r(\omega)\in(0,1)$. Then there exist a large enough $\delta=\delta(\omega,p)>0$ such that
\[\begin{split}\label{cdyl}
||D^{n}||^{q}_{A_{\omega}^{p}  \rightarrow L_{\mu}^{q}}\asymp\sup_{a\in\mathbb{D}}\int_{\mathbb{D}}\frac{(1-|a|)^{\delta q}}{ |1-\overline{a}w|^{(\delta+n)q} (\omega(S(a)))^{q/p} } d\mu(w).
\end{split}\]
\end{lem}

\begin{proof} For $a\in \mathbb{D}$ and $r\in(0,1)$,  we have
\[\begin{split}
\frac{\mu(\Delta(a,~r))}{(1-|a|)^{nq}(\omega(S(a)))^{q/p}}&=\int_{\Delta(a, r)}\frac{1}{ (1-|a|)^{nq} (\omega(S(a)))^{q/p}}d\mu(w)\\
&\asymp \int_{\Delta(a, r)} \frac{(1-|a|)^{\delta q}}{ |1-\overline{a}w|^{(\delta+n)q} (\omega(S(a)))^{q/p} }  d\mu(w)\\
&\lesssim \int_{\mathbb{D}}\frac{(1-|a|)^{\delta q}}{ |1-\overline{a}w|^{(\delta+n)q} (\omega(S(a)))^{q/p} } d\mu(w).\\
\end{split}\]
By Lemma \ref{D}, we find that
\[\begin{split}
||D^{n}||^{q}_{A_{\omega}^{p}  \rightarrow L_{\mu}^{q}}&\asymp\sup_{a\in\mathbb{D}}\frac{\mu(\Delta(a,~r))}{(1-|a|)^{nq}(\omega(S(a)))^{q/p}}\\
&\lesssim \sup_{a\in\mathbb{D}}\int_{\mathbb{D}}\frac{(1-|a|)^{\delta q}}{ |1-\overline{a}w|^{(\delta+n)q} (\omega(S(a)))^{q/p} } d\mu(w).\\
\end{split}\]
By \cite[Lemma 3.1]{P2016}, we can choose some large enough $\delta=\delta(\omega,p)>0$ and
\begin{equation}\label{cshs}
f_{a}(z)= \Big( \frac{1-|a|}{1-\overline{a}z}   \Big)^{\delta} \omega(S(a))^{-1/p},~~~~~~~~a, z\in \mathbb{D},
\end{equation}
then $||f_{a}||_{A_{\omega}^{p}}\lesssim 1$. Thus,  we get
\[\begin{split}
\int_{\mathbb{D}} |f^{(n)}_{a}(z)|^q d\mu(z)=\int_{\mathbb{D}}\frac{|a|^{n}(1-|a|)^{\delta q}}{ |1-\overline{a}z|^{(\delta+n)q} (\omega(S(a)))^{q/p} } d\mu(z)\lesssim ||D^{n}||^{q}_{A_{\omega}^{p}  \rightarrow L_{\mu}^{q}}.\\
\end{split}\]
The proof is complete.
\end{proof}

We use the pullback measure as an important tool to study the generalize weighted composition operators between different Bergman spaces with doubling weights. Let $\varphi$ be an analytic self-map of $\mathbb{D}$ and $0< q<\infty$. Assume that $u \in H(\mathbb{D})$, we define a finite positive Borel measure $\mu^\nu_{\varphi,~u}$ on $\mathbb{D}$ as follows:
\[
\mu^\nu_{\varphi,~u}(E)=\int_{\varphi^{-1}(E)}|u(z)|^q \nu(z)dA(z),
\]
where $E$ is a Borel subset of unit disk $\mathbb{D}$. For $D^{n}_{\varphi,u}: A^p_{\omega}\rightarrow A^q_{\nu}$, it can be clearly seen that
\begin{equation}\label{fhgs}
 ||D^{n}_{\varphi,u}f||_{A^q_{\nu}} =\int_{\mathbb{D}}|f^{(n)}(z)|^q d\mu^\nu_{\varphi,~u}(z),~~~~f\in A_{\omega}^p.
\end{equation}

%
\begin{lem}\label{sqzzhs}
Let $\omega \in \mathcal{D}$, $n\in \mathbb{N}\bigcup \{0\}$, $0<p<\infty$ and $0<r=r(\omega)<1$. If $f\in A^p_{\omega}$, there exists a constant $C=C(\omega)>0$ such that
$$|f^{(n)}(z)|^p\leq\frac{ C }{ \omega(S(z))}\int_{\Delta(z,r)}\frac{|f(w)|^p}{(1-|w|)^{np}} \widetilde{\omega}(w)dA(w)$$
for $z\in \mathbb{D}$. Here $\widetilde{\omega}(z)=\frac{\hat{\omega}(z)}{1-|z|}$.
\end{lem}

\begin{proof} It is clear that $  1-|z|\asymp 1-|w|$ for $w \in \Delta(z, r)$.  Since $\omega \in  \mathcal{D}$, by \cite[Lemma 2.1]{P2016}  and \cite[(2.27)]{PR2021}, there exist constants  $0<\alpha=\alpha(\omega)<\beta=\beta(\omega)<\infty$ and $C=C(\omega)\geq 1$ such that
\begin{equation}\label{2}
 \frac{1}{C}\left(\frac{1-r}{1-t}\right)^\alpha \leq \frac{\hat{\omega}(r)}{\hat{\omega}(t)}\leq C \left(\frac{1-r}{1-t}\right)^{\beta},
\end{equation}
where $0\leq r\leq t <1$. By (\ref{2}), we know that $\hat{\omega}(z)\asymp \hat{\omega}(w)$ for $w \in \Delta(z, r)$. 
By a direct calculation, we know that $\hat{\omega}(z)(1-|z|)\asymp \omega(S(z))$ for $\omega \in  \mathcal{D}$.
 By \cite[Lemma 2.1]{LK1985}, we claim that
\[\begin{split}
|f^{(n)}(z)|^p&\leq\frac{ C }{ (1-|z|)^{2+np}}\int_{\Delta(z,r)}|f(w)|^p dA(w)\\
&\asymp\frac{ C }{ \hat{\omega}(z)(1-|z|)(1-|z|)^{np}}\int_{\Delta(z,r)}|f(w)|^p \frac{\hat{\omega}(w)}{(1-|w|)}  dA(w)\\
&\asymp\frac{C}{\omega(S(z))}\int_{\Delta(z,r)}\frac{|f(w)|^p}{(1-|w|)^{np}} \widetilde{\omega}(w)dA(w).
\end{split}\]
\end{proof}

\begin{lem}\cite[Theorem 1.3]{L2021}\label{bdd2}
Let $0<p\leq q<\infty$ and $\omega,\nu \in  \mathcal{D}$. Assume that $\varphi$ is an analytic self-map of $\mathbb{D}$, $u \in A_{\nu}^{q} $ and $n\in \mathbb{N} \cup\{0\}$. Then $D^n_{\varphi,u} : A_{\omega}^{p}  \rightarrow A_{\nu}^{q}$ is bounded if and only if there exists a large enough $\delta=\delta(\omega,p)>0$ such that
\begin{equation}
\sup_{a\in \mathbb{D}} \int_{\mathbb{D}} \frac{(1-|a|)^{\delta q}|u(\xi)|^{q}\nu(\xi)}{|1-\overline{a}\varphi(\xi)|^{(\delta+n)q}(\omega(S(a)))^{q/p}}dA(\xi)<\infty.
\end{equation}

\end{lem}\msk

\section{order boundness of $D^n_{\varphi,u}: A^p_{\omega}\rightarrow A_{\nu}^q$ for $0<p,q<\infty$}

Next, we will study the order boundedness of $D^{n}_{\varphi, u} \colon ~~ A_{\omega}^{p} \to A_{\nu}^{q}$ for $0<p,q<\infty$.
\begin{thm}\label{U2}~~~~Let $0<p,~q<\infty$ and $\omega, \nu\in \mathcal{D}$. Suppose $n\in\mathbb{N}\cup \{0\}$. Let $\varphi$ be an analytic self-map of $\mathbb{D}$ and $u\in H(\D)$. Then $D^n_{\varphi,u}\colon ~~ A_{\omega}^{p} \to A_{\nu}^{q}$ is order bounded if and only if \\
\begin{equation}\label{xyjn}
  \int_\mathbb{D} \frac{|u(z)|^q \nu(z)}{(1-|\varphi(z)|^2)^{nq}(\omega(S(\varphi(z))))^{q/p} } dA(z)< \infty.
 \end{equation}
\end{thm}
\begin{proof}
Assume that $D^n_{\varphi,u}\colon ~~ A_{\omega}^{p} \to A_{\nu}^{q}$ is order bounded. There exists a non-negative function $h\in L^q_{\nu}$ such that $|D^n_{\varphi,u}f(z)|\leq h(z)$ for all $z\in \mathbb{D}$ and $f\in A_{\omega}^{p} $ with $\|f\|_{A_{\omega}^{p}}\lesssim1$.
To get (\ref{xyjn}), we set
\begin{equation}\label{orderI1}
\begin{split}
I_1=\int_{\{z\in\D,|\varphi(z)|>\frac{1}{2}\}}\frac{|u(z)|^q \nu(z)}{(1-|\varphi(z)|^2)^{nq}(\omega(S(\varphi(z))))^{q/p}} dA(z)\\
\end{split}
\end{equation}
and
\begin{equation}\label{orderI2}
\begin{split}
I_2=\int_{\{z\in\D,|\varphi(z)|\leq\frac{1}{2}\}}\frac{|u(z)|^q \nu(z)}{(1-|\varphi(z)|^2)^{nq}(\omega(S(\varphi(z))))^{q/p}} dA(z).\\
\end{split}
\end{equation}
By (\ref{cshs}) and for $z \in \mathbb{D}$, take\\
\begin{equation*}
f_{\varphi(z)}(w)=\frac{(1-|\varphi(z)|^2)^{\delta} }{ (1-\overline{\varphi(z)}w)^{\delta} (\omega(S(\varphi(z))))^{1/p} },~~~~w\in \mathbb{D}.
\end{equation*}
For some large enough $\delta=\delta(\omega,p)>0$,  we know that $f_{\varphi(z)}\in  A_{\omega}^{p}  $ and $\|f_{\varphi(z)}\|_{A_{\omega}^{p}}\lesssim 1$.
Then
\begin{equation}\label{GJDHS}
f_{\varphi(z)}^{(n)}(w)=C_{\delta, n}\frac{(1-|\varphi(z)|^2)^\delta (\overline{\varphi(z)})^n }{(1-\overline{\varphi(z)}w)^{\delta+n} (\omega(S(\varphi(z))))^{1/p}},
\end{equation}
where $C_{\delta, n}=\delta(\delta+1)(\delta+2)...(\delta+n-1)$.
By a direct computation,  for  $z \in \mathbb{D}$, we have
\[
|D^n_{\varphi,u}f_{\varphi(z)}(w)|=\frac{C_{\delta,n}(1-|\varphi(z)|^2)^\delta|u(w)||\varphi(z)|^n}{|1-\overline{\varphi(z)}\varphi(w) |^{n+\delta}  (\omega(S(\varphi(z))))^{1/p}}\leq h(w).
\]
So, by taking $w=z$, we can get
\[
\frac{C_{\delta,n}|u(z)||\varphi(z)|^n}{(1-|\varphi(z)|^2 )^{n}  (\omega(S(\varphi(z))))^{1/p}}= |D^n_{\varphi,u}f_{\varphi(z)}(z)|\leq h(z).
\]
For  $z\in \mathbb{D}$ such that $|\varphi(z)|>\frac{1}{2}$, we get $|\varphi(z)|^{n}>\frac{1}{2^n}$. Therefore,
\begin{equation}\label{order1}
\begin{split}
I_{1}&=\int_{\{z\in\D,|\varphi(z)|>\frac{1}{2}\}}\frac{|u(z)|^q}{(1-|\varphi(z)|^2)^{nq}(\omega(S(\varphi(z))))^{q/p}} \nu(z)dA(z)\\
&\leq\frac{2^{nq}}{C_{\delta,n}}\int_{_{\{z\in\D,|\varphi(z)|>\frac{1}{2}\}}}
\big|\frac{C_{\delta,n}|u(z)||\varphi(z)|^n}{(1-|\varphi(z)|^2)^{n}(\omega(S(\varphi(z))))^{1/p}}\big|^q \nu(z)dA(z)\\
&\lesssim\int_{\D}
\big|\frac{C_{\delta,n}|u(z)||\varphi(z)|^n}{(1-|\varphi(z)|^2)^{n}(\omega(S(\varphi(z))))^{1/p}}\big|^q \nu(z)dA(z)\\
&\leq\int_{\D}|h(z)|^q \nu(z)dA(z)<\infty.
\end{split}
\end{equation}
For $z\in \mathbb{D}$ such that $|\varphi(z)|\leq \frac{1}{2}$, 
we can find a constant $C>0$ such that
\begin{equation}\label{AA2}
\frac{1}{(1-|\varphi(z)|^2)^{n} (\omega(S(\varphi(z))))^{1/p}}\leq C.
\end{equation}
On the other hand, since $P_{n}(z)=\frac{z^n}{\|z^n\|_{A^p_{\omega}}}$ is in $A^p_{\omega}$ and $\|P_{n}\|_{A^p_{\omega}}\leq 1$, by the order boundedness of the operator $D^n_{\varphi,u}$, for $z\in \mathbb{D}$, we obtain\\
\begin{equation}\label{AA}
\frac{n!}{\|z^n\|_{A^p_{\omega}}}|u(z)|= |D^n_{\varphi,u}P_{n}(z) | \leq h(z).
\end{equation}
Since $n$ is fixed, from (\ref{AA2}) and (\ref{AA}), for $z\in \mathbb{D}$, we get
\begin{equation}\label{order2}
\begin{split}
I_{2}&=\int_{\{z\in\D,|\varphi(z)|\leq \frac{1}{2}\}} \frac{|u(z)|^q}{(1-|\varphi(z)|^2)^{nq} (\omega(S(\varphi(z))))^{q/p}}\nu(z)dA(z)\\
&\leq C\int_{\{z\in\D,|\varphi(z)|\leq \frac{1}{2}\}}|u(z)|^q \nu(z)dA(z)\lesssim\int_{\D}|u(z)|^q \nu(z)dA(z)\\
&\leq\int_{\D}|h(z)|^q \nu(z)dA(z)<\infty.
\end{split}
\end{equation}

By (\ref{order1}) and (\ref{order2}), we see that
$$  \int_{\mathbb{D}}\frac{|u(z)|^q}{(1-|\varphi(z)|^{2})^{nq} (\omega(S(\varphi(z))))^{q/p}} \nu(z)dA(z)=I_{1}+I_{2}<\infty.$$
Thus, the condition (\ref{xyjn}) holds.

Conversely, assume that condition (\ref{xyjn}) holds. Define
$$h(z)=\frac{|u(z)|}{(1-|\varphi(z)|^2)^{n} (\omega(S(\varphi(z))))^{1/p}}.$$
 Then $h$ is a nonnegative function in $L^q_{\nu}$. For any function $f\in A^p_{\omega}$ with $||f||_{A^p_{\omega}}\leq 1$, by Lemma \ref{DZ}, there is a constant $C=C(\omega)>0$ such that
\[
|D^n_{\varphi,u}f(z)|=|u(z)f^{(n)}(\varphi(z))|\leq C \frac{|u(z)|}{(1-|\varphi(z)|^2)^{n} (\omega(S(\varphi(z))))^{1/p}}=Ch(z)
\]
for any $z\in \mathbb{D}$. Thus, $ D^n_{\varphi,u}: A^p_{\omega} \rightarrow A^q_{\nu}$ is order bounded. The proof is complete.
\end{proof}


\section{Essential norm of $D^n_{\varphi,u}: A^p_{\omega}\rightarrow A_{\nu}^q$ for $1\leq p\leq q<\infty$}
We begin this section with an approximation of the essential norm of the bounded operator $D^n_{\varphi,u}: A^p_{\omega}\rightarrow A_{\nu}^q$ for $1\leq p\leq q<\infty$. If $f\in H(\mathbb{D})$, then $f(z)= \sum_{k=0}^{\infty}a_{k}z^k$. For any $m\geq 1$, let $R_{m}f(z)= \sum_{k=m}^{\infty}a_{k}z^k$  and $T_{m}=I-R_{m}$, where $If = f $ is the identity operator. In order to prove one of the main results, we need the following lemmas.
\begin{lem}\cite[Proposition 1]{zk1991}\label{mt}
Suppose $X$ is a Banach space of holomorphic functions in $\mathbb{D}$ with the property that the polynomials are dense in $X$. Then $||T_{m}f-f||_{X}\rightarrow 0$ as $m\rightarrow \infty$ for each $f\in X$ if and only if $\sup\{||T_{m}||: m\geq1\}<\infty$.
\end{lem}

\begin{lem}\cite[Corollary 3]{zk1991}\label{lemm1}
The Taylor series of every function in $H^p$ converges in norm if and only if $1<p<\infty$.
\end{lem}

\begin{lem}\label{bxyl}
For $1< p<\infty$ and $\omega$ is a radial weight, then $\|T_{m}f-f\|_{A^p_{\omega}}\rightarrow 0$ as $m\rightarrow \infty$ for each $f\in A^p_{\omega}$. Moreover, $\sup\{\|R_{m}\|_{A^{p}_{\omega} \rightarrow A^{p}_{\omega} }: m\geq1\}<\infty$ and $\sup\{\|T_{m}\|_{A^{p}_{\omega} \rightarrow A^{p}_{\omega} }: m\geq1\}<\infty$, where $R_{m}=I-T_{m}$.
\end{lem}

\begin{proof}
It follows from Lemmas \ref{mt} and  \ref{lemm1} that $T_{m}$ is bounded uniformly on $H^p$ for $1<p<\infty$.
Thus, there exists a constant $C>0$ such that
\[
\frac{1}{2\pi}\int_{0}^{2\pi}|T_{m}f(re^{i\theta}) |^p d\theta \leq C \frac{1}{2\pi}\int_{0}^{2\pi}|f(re^{i\theta}) |^p d\theta,
\]
for $p>1$ and any $m\geq 1$. Applying polar coordinates, we see that
\[
||T_{m}f||_{A^p_{\omega}}^p \leq C \int_{0}^{1}\omega(r)r dr\int_{0}^{2\pi}|f(re^{i\theta}) |^p d\theta\leq C||f||_{A_{\omega}^p}^p.
\]
Therefore $||T_{m}||_{A^{p}_{\omega} \rightarrow A^{p}_{\omega} }\leq C$ for any $m\geq1$. By Lemma \ref{mt}, we obtain that $||T_{m}f-f||_{A_{\omega}^p}\rightarrow 0$ as $m\rightarrow \infty$. Since $R_{m}=I-T_{m}$, we have
\[
||R_{m}||_{A^{p}_{\omega} \rightarrow A^{p}_{\omega} }= ||I-T_{m}||_{A^{p}_{\omega} \rightarrow A^{p}_{\omega}}\leq 1+||T_{m}||_{A^{p}_{\omega} \rightarrow A^{p}_{\omega} }\leq 1+C.
\]
\end{proof}

\begin{lem}\label{p>1}
Suppose that $\omega \in \hat{\mathcal{D}}$ and $1<p<\infty$. Let $\varepsilon >0$ and $r\in(0,1)$. Then there exists a $m_{0}\in \mathbb{N}$, for any $m\geq m_{0},$
\begin{equation}
 |R_{m}f(z)|\lesssim\varepsilon ||f||_{A^p_{\omega}},
\end{equation}
for every $z\in D_{r}=\{z\in \mathbb{D}, |z|\leq r\}$ and each $f\in A^p_{\omega}$.
\end{lem}

\begin{proof}
Let $\omega_{n}=\int_{0}^1 r^n \omega(r)dr$. By \cite[P665]{PRS2018}, we see that $$B_{z}^{\omega}(\xi)=\sum_{n=0}^{\infty}\frac{(\xi\overline{z})^n}{2\omega_{2n+1}}$$ is the reproducing kernel of $A^p_{\omega}$ for  $p\geq1$. 
Then, we have
\begin{equation}\label{3}
\begin{split}
|R_{m}f(z)|=|\langle R_{m}f, B_{z}^{\omega}\rangle|
=|\langle f, R_{m}B_{z}^{\omega}\rangle|&\leq\int_{\mathbb{D}}|f(w)\overline{R_{m}B_{z}^{\omega}(w) }| \omega(w)dA(w)\\
&\lesssim ||f||_{A^p_{\omega}}|| R_{m}B_{z}^{\omega}||_{A^q_{\omega}},
\end{split}
\end{equation}
where $\frac{1}{p}+\frac{1}{q}=1$. For $z\in D_{r}$, we show that
\[
|| R_{m}B_{z}^{\omega}||_{A^q_{\omega}} = \left(\int_{\mathbb{D}}|R_{m}B_{z}^{\omega}(w) |^q \omega(w)dA(w)\right)^{\frac{1}{q}}\lesssim \sum_{k=m}^{\infty}\frac{r^{k}}{2\omega_{2k+1}}.
\]
By \cite[Lemma 6]{PRS2018}, we deduce that
\begin{equation}\label{qh}
\lim_{m\rightarrow\infty} \sum_{k=m}^{\infty}\frac{r^{k}}{2\omega_{2k+1}}=0.
\end{equation}
 Therefore, for any $\varepsilon>0$, there exists a $m_{0}\in \mathbb{N}$ and $m\geq m_{0}$, such that
 $$|| R_{m}B_{z}^{\omega}||_{A^p_{\omega} } \leq  \varepsilon. $$
By (\ref{3}), we get $|R_{m}f(z)|\lesssim\varepsilon ||f||_{A^p_{\omega}}$ for any $f\in A^p_{\omega}$.
\end{proof}

For $p=1$, let $\mathcal{T}_{m}f(z)=\sum_{k=0}^{m-1}(1-\frac{k}{m}) a_{k} z^k $ and $\mathcal{R}_{m}=I-\mathcal{T}_{m}$. We get the following lemma.

\begin{lem}\label{lp=1}
Let $\omega \in \hat{\mathcal{D}}$ and $f\in A^1_{\omega}$, then $||\mathcal{T}_{m}f-f||_{A^1_{\omega}}\rightarrow 0$ as $m\rightarrow \infty$ for each $f\in A^1_{\omega}$. Moreover,
  $\sup\{||\mathcal{T}_{m}||_{A^1_{\omega}\rightarrow A^1_{\omega}}: m\geq1\}<\infty$  and $\sup\{||\mathcal{R}_{m}||_{A^1_{\omega}\rightarrow A^1_{\omega}}: m\geq1\}<\infty$, where $\mathcal{R}_{m}=I-\mathcal{T}_{m}$.
\end{lem}

\begin{proof}
  By \cite[P.196]{D2011},   $||\mathcal{T}_{m}||_{H^1 \rightarrow H^1 }\leq1$.  Using the same way of Lemma \ref{bxyl},
we know that $||\mathcal{T}_{m}||_{A^{1}_{\omega} \rightarrow A^{1}_{\omega} }\leq C$ for any $m\geq1$.
We claim that 
\begin{equation}\label{p1-note}
||\mathcal{R}_{m}||_{A^{1}_{\omega} \rightarrow A^{1}_{\omega} }= ||I-\mathcal{T}_{m}||_{A^{1}_{\omega} \rightarrow A^{1}_{\omega}}\leq 1+||\mathcal{T}_{m}||_{A^{1}_{\omega} \rightarrow A^{1}_{\omega} }< 1+C
\end{equation}
for any $m\geq1$ and C is a positive constant.
\end{proof}

\begin{lem}\label{p=1}
Assume that $\omega \in \hat{\mathcal{D}}$. Let $\varepsilon >0$ and $r\in(0,1)$. Then there exists a $m_{0}\in \mathbb{N}$, for any $m\geq m_{0},$
\begin{equation}
 |\mathcal{R}_{m}f(w)|\lesssim \varepsilon ||f||_{A^1_{\omega}},
\end{equation}
for every $w\in D_{r}=\{w\in \mathbb{D}, |w|\leq r\}$ and each $f\in A^1_{\omega}$.
\end{lem}

\begin{proof}
By the proof of Lemma \ref{p>1},
we deduce that
\[
|\mathcal{R}_{m}f(w)|=|\langle \mathcal{R}_{m}f, B_{w}^{\omega}\rangle|=|\langle f, \mathcal{R}_{m}B_{w}^{\omega}\rangle|
\lesssim ||f||_{A^1_{\omega}}|| \mathcal{R}_{m}B_{w}^{\omega}||_{H^{\infty}}.
\]
Take $|w|\leq r$, we can prove that
\[
|| \mathcal{R}_{m}B_{w}^{\omega}||_{H^{\infty}} = \sup_{\xi\in \mathbb{D}}|\mathcal{R}_{m}B_{w}^{\omega}(\xi) | = \sup_{\xi\in \mathbb{D}}|(I-\mathcal{T}_{m})B_{w}^{\omega}(\xi) |\leq \frac{1}{m}\sum_{k=1}^{\infty}\frac{kr^{k-1}}{2\omega_{2k+1}} + \sum_{k=m}^{\infty}\frac{r^{k}}{2\omega_{2k+1}}.
\]
By \cite[Lemma 6]{PRS2018}, we see that $\sum_{k=1}^{\infty}\frac{kr^{k-1}}{2\omega_{2k+1}}$ is convergent and (\ref{qh}) holds. 
Therefore, for any $\varepsilon>0$, there exists a $m_{0}\in \mathbb{N}$ and $m\geq m_{0}$, such that
 $$|| \mathcal{R}_{m}B_{w}^{\omega}||_{H^{\infty}} \leq  \varepsilon. $$
Thus $|\mathcal{R}_{m}f(w)|\lesssim \varepsilon ||f||_{A^1_{\omega}}$ for any $f\in A^1_{\omega}$.
\end{proof}

The following lemma is very useful to prove the compactness of composition operators and its generalizations on some function spaces.  

\begin{lem}\cite[Lemma 2.2]{L2021}\label{cpt}
Suppose $0< p,\,q<\infty$, $\omega, \nu \in  \mathcal{D}$. Suppose $u\in H(\mathbb{D})$ and $n\in \mathbb{N}\bigcup\{0\}$. Let $\varphi$ be an analytic self-map of $\mathbb{D}$ such that  $ D^{n}_{\varphi,u}: A^{p}_{\omega} \rightarrow A^{q}_{\nu}$ is bounded. Then $D^{n}_{\varphi,u}: A^{p}_{\omega} \rightarrow A^{q}_{\nu}$ is compact if and only if whenever $\{f_k\}$ is
bounded in $A_{\omega}^{p}$ and $f_k\rightarrow 0$ uniformly on compact subsets of
$\mathbb{D}$ as $k\rightarrow \infty$, $\lim_{k\rightarrow \infty}\|D^{n}_{\varphi,u}(f_k)\|_{A^{q}_{\nu}}=0$.
\end{lem}

\begin{thm}\label{ess}
Let $1\leq p\leq q<\infty$ and $\omega, \nu \in  \mathcal{D}$. Suppose  $n\in \mathbb{N}\bigcup\{0\}$. Let $\varphi$ be an analytic self-map of $\mathbb{D}$ and $u \in A^{q}_{\nu}$. If $D^{n}_{\varphi,u} : A_{\omega}^{p}  \rightarrow A^{q}_{\nu}$ is bounded, then there exists a large enough $\delta=\delta(\omega,p)>0$ such that
\begin{equation}\label{jl}
||D^{n}_{\varphi,u}||_{e,~A^p_{\omega}\rightarrow A^{q}_{\nu}}^{q} \asymp\limsup_{|a|\rightarrow 1} \int_{\mathbb{D}} \frac{(1-|a|)^{\delta q}|u(\xi)|^{q} \nu(\xi) }{|1-\overline{a}\varphi(\xi)|^{(\delta+n)q}(\omega(S(a)))^{q/p}}dA(\xi).
\end{equation}
\end{thm}

\begin{proof}
\textbf{Lower estimate}.
Let $f_{a}(z)= \Big( \frac{1-|a|}{1-\overline{a}z}   \Big)^{\delta} \omega(S(a))^{-1/p}$ for some large enough $\delta=\delta(\omega,p)>0$. Then $\{f_{a}\}$ is a bounded sequence in $A^p_{\omega}$ converging to zero uniformly on compact subsets of $\mathbb{D}$ as $ |a| \rightarrow 1$.
Fix a compact operator $K:~~A^p_{\omega} \rightarrow A^{q}_{\nu}$, by Lemma \ref{cpt}, we know that $||Kf_{a}||_{A^{q}_{\nu}}\rightarrow 0$ as $ |a| \rightarrow 1$. Therefore
\[\begin{split}
||D^{n}_{\varphi,u}-K||_{A^p_{\omega} \rightarrow A^{q}_{\nu} } &\gtrsim \limsup_{|a|\rightarrow 1}||(D^{n}_{\varphi,u}-K)f_{a}||_{A^{q}_{\nu}}\\
&\gtrsim \limsup_{|a|\rightarrow 1}(||D^{n}_{\varphi,u} f_{a}||_{A^{q}_{\nu}}-||Kf_{a}||_{A^{q}_{\nu}} ) \\
&= \limsup_{|a|\rightarrow 1}||D^{n}_{\varphi,u} f_{a}||_{A^{q}_{\nu}}.
\end{split}\]
Moreover, we have
 \[\begin{split}
 ||D^{n}_{\varphi,u}||_{e,~A^p_{\omega} \rightarrow A^q_{\nu}} &= \inf_{K}||D^{n}_{\varphi,u}-K||_{A^p_{\omega} \rightarrow A^{q}_{\nu}}\gtrsim \limsup_{|a|\rightarrow 1}||D^{n}_{\varphi,u} f_{a}||_{A^{q}_{\nu}}\\
 &=\limsup_{|a|\rightarrow 1} 
 \left(\int_{\mathbb{D}} \frac{|a|^{n}(1-|a|)^{\delta q}|u(\xi)|^{q} \nu(\xi) }{|1-\overline{a}\varphi(\xi)|^{(\delta+n)q}(\omega(S(a)))^{q/p}}dA(\xi)\right)^{1/q}.
 \end{split}\]
We get
\begin{align}\label{blew}
||D^{n}_{\varphi,u}||^q_{e,~A^p_{\omega} \rightarrow A^{q}_{\nu}}&\gtrsim\limsup_{|a|\rightarrow 1} \int_{\mathbb{D}} \frac{(1-|a|)^{\delta q}|u(\xi)|^{q} \nu(\xi) }{|1-\overline{a}\varphi(\xi)|^{(\delta+n)q}(\omega(S(a)))^{q/p}}dA(\xi).
\end{align}

\textbf{Upper estimate}. The case $1< p\leq q <\infty$.
Considering the compact operator $T_{m}: A_{\omega}^{p} \rightarrow A^{q}_{\nu}$ by $T_{m}f=\sum_{k=0}^{m-1}b_{k}z^{k}$
and letting $R_{m}= I - T_{m}$,  where $I$ is identity operator.
We can see that
\[
||D^{n}_{\varphi,u}||_{e,~A^p_{\omega} \rightarrow A^{q}_{\nu}} \leq  ||D^{n}_{\varphi,u}\circ R_{m} ||_{e,~A^p_{\omega} \rightarrow A^{q}_{\nu}} + ||D^{n}_{\varphi,u}\circ   T_{m}||_{e,~A^p_{\omega} \rightarrow A^{q}_{\nu}} =  ||D^{n}_{\varphi,u}\circ R_{m} ||_{e,~A^p_{\omega} \rightarrow A^{q}_{\nu}} .
\]
Thus
\begin{equation}\label{syj}
||D^{n}_{\varphi,u}||^q_{e,~A^p_{\omega} \rightarrow A^{q}_{\nu}}\leq \liminf_{m\rightarrow \infty}||D^{n}_{\varphi,u}\circ R_{m}||^q_{e,~A^p_{\omega} \rightarrow  A^{q}_{\nu}}
 \leq\liminf_{m\rightarrow \infty} ||D^{n}_{\varphi,u}\circ R_{m} ||^q_{A^p_{\omega} \rightarrow A^{q}_{\nu}}.
\end{equation}
Fix $f\in A_{\omega}^{p} $ with $||f||_{A^p_\omega}\leq 1$ and $r \in (0,1)$. Suppose $D_{r}=\{z\in \mathbb{D},~|z|\leq r\}$. Then
 \[\begin{split}
||(D^{n}_{\varphi,u}\circ R_{m})f||^q_{A^p_{\omega} \rightarrow A_{\nu}^q} &\leq \int_{\mathbb{D}}|(R_{m}f)^{(n)}(\varphi(\xi))|^{q} |u(\xi)|^{q} \nu(\xi)dA(\xi)\\
&=\int_{\mathbb{D}}|(R_{m}f)^{(n)}(z)|^{q}d\mu^\nu_{\varphi,~u}(z),\\
\end{split}\]
where $\mu^\nu_{\varphi,~u}=\int_{\varphi^{-1}(E)}|u(z)|^q \nu(z)dA(z)$ for all $E$ is Borel subsets of $\mathbb{D}$.

From Lemma \ref{DZ}, we have
\begin{equation}\label{q-p}
   |f(z)|^{q-p}
     \lesssim \frac{||f||^{q-p}_{A_{\omega}^p}}{(\omega(S(z)))^{(q-p)/p}}.\\
  \end{equation}
By Lemma \ref{sqzzhs} and (\ref{q-p}),  we obtain
\begin{equation}\label{yyds}
\begin{split}
&\int_{\mathbb{D}}|(R_{m}f)^{(n)}(z)|^{q}d\mu^\nu_{\varphi,~u}(z)\\
\leq& \int_{\mathbb{D}} d\mu^\nu_{\varphi,~u}(z)\frac{C}{\omega(S(z))}\int_{\Delta(z,r)}\frac{|R_{m}f(w)|^{q} }{(1-|w|)^{nq}}\widetilde{\omega}(w) dA(w)\\
\asymp &  \int_{\mathbb{D}} d\mu^\nu_{\varphi,~u}(z)\int_{\Delta(z,r)}\frac{|R_{m}f(w)|^{q-p+p}}{(1-|w|)^{nq} \omega(S(w))}\widetilde{\omega}(w)dA(w)\\
\lesssim& ||R_{m}f||^{q-p}_{A^{p}_{\omega}}\int_{\mathbb{D}} d\mu^\nu_{\varphi,~u}(z)\int_{\Delta(z,r)}\frac{|R_{m}f(w)|^{p}}{ (1-|w|)^{nq}(\omega(S(w)))^{(q-p)/p} \omega(S(w))}\widetilde{\omega}(w)dA(w)\\
=& ||R_{m}f||^{q-p}_{A^{p}_{\omega}}\int_{\mathbb{D}} d\mu^\nu_{\varphi,~u}(z)\int_{\Delta(z,r)}\frac{|R_{m}f(w)|^{p}}{ (1-|w|)^{nq}(\omega(S(w)))^{q/p} } \widetilde{\omega}(w) dA(w)\\
=&  ||R_{m}f||^{q-p}_{A^{p}_{\omega}}\int_{\mathbb{D}} d\mu^\nu_{\varphi,~u}(z)\int_{\mathbb{D}}\frac{\chi_{\Delta(z,r)}(w)|R_{m}f(w)|^{p}}{ (1-|w|)^{nq}(\omega(S(w)))^{q/p} } \widetilde{\omega}(w) dA(w),
\end{split} \end{equation}
where $\chi_{\Delta(z,r)}$ is the characteristic function of the set $\Delta(z, r)$. Obviously, $\chi_{\Delta(z,r)}(w)=\chi_{\Delta(w,r)}(z)$. By Fubini's Theorem,   
we obtain
\begin{equation}\label{yygs}
\begin{split}
&\int_{\mathbb{D}}|(R_{m}f)^{(n)}(z)|^q d\mu^\nu_{\varphi,~u}(z)\\ \lesssim& ||R_{m}f||^{q-p}_{A^{p}_{\omega}}\int_{\mathbb{D}} \frac{\mu^\nu_{\varphi,~u}(\Delta(w,r))}{ (1-|w|)^{nq}(\omega(S(w)))^{q/p} } |R_{m}f(w)|^{p}\widetilde{\omega}(w) dA(w).\\
\end{split}\end{equation}
Set
\[\begin{split}
J_{1,~m}&=\int_{D_r}\frac{\mu^\nu_{\varphi,~u}(\Delta(w,r))}{ (1-|w|)^{nq}(\omega(S(w)))^{q/p} }|R_{m}f(w)|^p\widetilde{\omega}(w) dA(\xi)\\
\end{split}\]
and
\[\begin{split}
J_{2,~m}&=\int_{\mathbb{D}\backslash D_{r}}
\frac{\mu^\nu_{\varphi,~u}(\Delta(w,r))}{(1-|w|)^{nq}  (\omega(S(w)))^{q/p}   }|R_{m}f(w)|^p\widetilde{\omega}(w) dA(w).\\
\end{split}\]
Then we get
\begin{equation}\label{jh}
||(D^{n}_{\varphi,u}\circ R_{m})f||^q_{A^p_{\omega} \rightarrow A^q_{\nu} }\lesssim||R_{m}f||^{q-p}_{A_{\omega}^p} (  J_{1,~m} +  J_{2,~m} ),
\end{equation}
with $m\geq 1$. Since $D^{n}_{\varphi,u}: A^p_\omega \rightarrow A^q_{\nu}$ is bounded, Lemma \ref{bdd2} implies that there exists a large enough $\delta=\delta(\omega,p)>0$ such that
\begin{align}\label{MM}
M&=\sup_{a\in\mathbb{D}}\int_{\mathbb{D}} \frac{(1-|a|)^{\delta q} |u(w)|^{q} \nu(w) }{|1-\overline{a}\varphi(w)|^{(\delta+n)q}(\omega(S(a)))^{q/p}}dA(w)\nonumber\\
&=\sup_{a\in\mathbb{D}}\int_{\mathbb{D}} \frac{(1-|a|)^{\delta q} }{|1-\overline{a}\xi|^{(\delta+n)q}(\omega(S(a)))^{q/p}}d\mu^\nu_{\varphi,~u}(\xi)\nonumber\\
&< \infty.
\end{align}
For $\xi\in \Delta(w,r)$, we have
\begin{align}\label{yl2.4}
\frac{\mu^\nu_{\varphi,~u}(\Delta(w,~r))}{(1-|w|)^{nq}(\omega(S(w)))^{q/p}}&= \int_{\Delta(w, r)}\frac{1}{ (1-|w|)^{nq} (\omega(S(w)))^{q/p}}d\mu^\nu_{\varphi,~u}(\xi)\nonumber\\
&\asymp \int_{\Delta(w, r)} \frac{(1-|w|)^{\delta q}}{ |1-\overline{w}\xi|^{(\delta+n)q} (\omega(S(w)))^{q/p} }  d\mu^\nu_{\varphi,~u}(\xi)\nonumber\\
&\lesssim \int_{\mathbb{D}}\frac{(1-|w|)^{\delta q}}{ |1-\overline{w}\xi|^{(\delta+n)q} (\omega(S(w)))^{q/p} } d\mu^\nu_{\varphi,~u}(\xi).\nonumber\\
\end{align}
Fix $\varepsilon > 0$.
By (\ref{yl2.4}), (\ref{MM}) and Lemma \ref{p>1}, hence
\[\begin{split}
J_{1,~m}&=\int_{D_r}\frac{\mu^\nu_{\varphi,~u}(\Delta(w,r))}{(1-|w|)^{nq}  (\omega(S(w)))^{q/p}   }|R_{m}f(w)|^p \widetilde{\omega}(w)dA(w)\\
&\leq\sup_{w\in\mathbb{ D}}\frac{\mu^\nu_{\varphi,~u}(\Delta(w,r))}{(1-|w|)^{nq}  (\omega(S(w)))^{q/p}   } \int_{D_r} |R_{m}f(w)|^p \widetilde{\omega}(w)dA(w)\\
&\leq CM\int_{D_r}|R_{m}f(w)|^p\widetilde{\omega}(w)dA(w)\\
&\leq CM \varepsilon^p ||f||^p_{A^p_{\omega}},
\end{split}\]
for any $m\geq m_{0}$. 
Thus,
\begin{equation}\label{1n}
\lim_{m\rightarrow\infty} \sup_{||f||_{A^p_\omega}\leq 1}||R_{m}f||^{q-p}_{A_{\omega}^p}  J_{1,~m}=0.
\end{equation}
For $\omega \in \mathcal{D}$ and $f\in H(\mathbb{D})$, from \cite[Proposition 5]{PRS2018}, we know that
  \begin{equation}\label{1}
 ||f||_{A^{p}_{\widetilde{\omega}}} \asymp ||f||_{A^{p}_{\omega}}.
 \end{equation}
By (\ref{yl2.4}), (\ref{1}) and Lemma \ref{bxyl}, 
we claim that
\[\begin{split}
J_{2,~m}&=\int_{\mathbb{D}\backslash D_{r}}
\frac{\mu^\nu_{\varphi,~u}(\Delta(w,r))}{(1-|w|)^{nq}  (\omega(S(w)))^{q/p}   }|R_{m}f(w)|^p\widetilde{\omega}(w) dA(w)\\
&\leq \sup_{|a|>r }
\frac{\mu^\nu_{\varphi,~u}(\Delta(a,r))}{(1-|a|)^{nq}  (\omega(S(a)))^{q/p}}\int_{\mathbb{D}\backslash D_{r}}|R_{m}f|^{p}\widetilde{\omega}(w)dA(w)\\
&\lesssim \sup_{m\geq1}||R_{m}||_{A^{p}_{\omega} \rightarrow A^{p}_{\omega} }^p ||f||^p_{A^p_\omega} \sup_{|a|>r}\int_{\mathbb{D}} \frac{(1-|a|)^{\delta q}}{|1-\overline{a}\xi|^{(\delta+n)q}(\omega(S(a)))^{q/p}}d\mu^\nu_{\varphi,~u}(\xi).
\end{split}\]
Hence,
\begin{equation}\label{2n}
\lim_{m\rightarrow\infty} \sup_{||f||_{A^p_\omega}\leq 1}||R_{m}f||^{q-p}_{A_{\omega}^p}  J_{2,~m}\lesssim \sup_{|a|>r}\int_{\mathbb{D}} \frac{(1-|a|)^{\delta q}}{|1-\overline{a}\xi|^{(\delta+n)q}(\omega(S(a)))^{q/p}}d\mu^\nu_{\varphi,~u}(\xi).
\end{equation}
Combining (\ref{syj}), (\ref{jh}), (\ref{1n}), (\ref{2n}) and (\ref{MM}), we deduce that
\[
||D^{n}_{\varphi,u}||_{e,~A^p_{\omega} \rightarrow  A^q_{\nu} }^q
 \leq\liminf_{m\rightarrow \infty} ||D^{n}_{\varphi,u}\circ R_{m} ||^q_{A^p_{\omega} \rightarrow A_{\nu}^q}\]
 \[\lesssim\sup_{|a|>r}\int_{\mathbb{D}} \frac{(1-|a|)^{\delta q}|u(w)|^{q} \nu(w) }{|1-\overline{a}\varphi(w)|^{(\delta+n)q}\omega(S(a))^{q/p}}dA(w).
\]
Letting $r\rightarrow 1$, we have
\begin{equation}\label{upp}
||D^{n}_{\varphi,u}||_{e,~A^p_{\omega} \rightarrow A^q_{\nu}}^q \lesssim \limsup_{|a|\rightarrow 1} \int_{\mathbb{D}} \frac{(1-|a|)^{\delta q}|u(w)|^{q}\nu(w) }{|1-\overline{a}\varphi(w)|^{(\delta+n)q}\omega(S(a))^{q/p}}dA(w).
\end{equation}

When $1=p\leq q <\infty$, by Lemmas \ref{lp=1} and   \ref{p=1},  we can use the same way to get that   (\ref{upp}) holds. We omit the details.
The proof of the Theorem \ref{ess} is complete.
\end{proof}


\vspace{0.3truecm}
\noindent {\bf Data Availability. }

All data generated or analyzed during this study are included in this article and in its bibliography.

\vspace{0.1truecm}
\noindent {\bf Conflict of Interest. }

The authors declared that they have no conflict of interest.

\noindent {\bf Acknowledgements. }
\thanks{The author is extremely thankful to Professor Hasi Wulan  for his kind suggestions. The research was supported by National Natural Science Foundation of China (Nos.11720101003 and 12171299) and  Guangdong Basic and Applied Basic Research Foundation (No.2022A1515012117). The author would like to thank the anonymous referee for his/her careful reading of the manuscript and  valuable comments.}

\end{document}